\documentclass[12pt]{article}
\title{ On $\psi$- basic Bernoulli-Ward polynomials}
\author{A.K.Kwa\'sniewski \\
\\Higher School of Mathematicand Applied Informatics\\
PL-15-021 Bia{\l}ystok, ul.Kamienna 17, POLAND\\
e-mail: kwandr@uwb.edu.pl}

\usepackage{amsmath,amsthm}

\chardef\bslash=`\\ 
\begin{document}
\maketitle
\begin{abstract}
The  Ward solution  of $\psi$- {\em difference calculus}
nonhomogeneous equation
$$\Delta_{\psi}f=\varphi\;\;\;\;\;\;\varphi=?$$
is found in the form of
$$f(x)=\sum_{n \geq
1}\frac{B_{n}}{n_{\psi}!}\varphi^{(n-1)}(x)+\int_{\psi}\varphi(x)
+p(x)$$ (where $B_{n}$ denote $\psi$-{\em Bernoulli-Ward numbers}
\cite{1}) - in the framework of the $\psi$-{\em Finite Operator
Calculus} \cite{2} - \cite{5}. Specifications to $q$-calculus case
and the new Fibonomial calculus case \cite{5,6} are made explicit.
\end{abstract}
A.M.S Classification numbers: 11B39, 11B65, 05A15
\section{Remark on Notation and References}
At first let us anticipate with $\psi$- remark. $\psi$ denotes an
extension of
$$\left\{\frac{1}{n!}\right\}_{n \geq 0}$$
sequence to quite arbitrary one ("admissible") and the specific
choices are for example: Fibonomially - extended ($F_{n},\;\;n\geq
0$ - Fibonacci sequence) Gauss $q$- extended
$$ \left\{\psi_{n}\right\}_{n \geq 0}=\left\{\frac{1}{F_{n}!}\right\}_{n \geq 0},\;\;\;\;
\left\{\psi_{n}\right\}_{n \geq
0}=\left\{\frac{1}{n_{q}!}\right\}_{n \geq 0}$$ admissible
sequences of extended umbral operator calculus - see more below.
With such an extension we may $\psi$ - mnemonic repeat with
exactly the same simplicity and beauty much of what was done by
Rota years ago. Thus due to efficient usage we get used to write
down these extensions in mnemonic upside down notation \cite{2,5}
$$n_{\psi}\equiv \psi_{n},\;\;\; x_{\psi}\equiv \psi (x) \equiv
\psi_{x},\;\;\;n_{\psi}!=n_{\psi}(n-1)_{\psi}!,\;\;\;0_{\psi}=1$$
$$x_{\psi}^{\underline{k}}=x_{\psi}(x-1)_{\psi}
\ldots(x-k+1)_{\psi}\equiv \psi (x)\psi (x-1) \ldots \psi(x-k+1)$$
You may consult for further development and use of this notation
\cite{4,5} and references therein.Summing up - we say it again.
The papers of main reference are \cite{1,2,3}. For not only
mnemonic reasons we follow here the notation from \cite{2,3} and
we shall take the results from \cite{1} as well as from \cite{2,3}
- for granted. (Note the access also via ArXiv to \cite{3,5}). For
other respective references see: \cite{2,3,5}. Note that we use
$\psi$-extension symbols - as in popular $q$-calculi (for example
$n_{q}$) - in upside down notation; for example:
$${\bf \psi_{n}\equiv n_{q}\; , \;\;\; x_{q}=\psi(x)=\frac{1-q^{x}}{1-q}.}$$
$B_{n}$ denote here $\psi$-{\em Bernoulli-Ward numbers}.
\section{}
Let the $\psi$-{\em difference delta operator} be defined as
$\Delta_{\psi}=E^{y}(\partial_{\psi})-id$. Then $\psi$- {\em basic
Bernoulli-Ward polynomials} $\left\{ B_{n}(x)\right\}_{n\geq 0}$
might be defined equivalently by
\begin{equation}\label{psi_B_1}
\sum_{s=0}^{n-1}\binom{n}{s}_{\psi}B_{s}(x)=nx^{n-1}\; ; n\geq 2\;
;\;\; B_{0}(x)=1,
\end{equation}
where $B_{n}=B_{n}(0)$ denote $\psi$-{\em Bernoulli-Ward numbers}:
$\left\{ B_{s}\right\}_{s \geq 0}$ or via
\begin{equation}\label{psi-B-2}
B_{n}(x)=\sum_{s=0}^{n}\binom{n}{s}_{\psi} B_{s}x^{n-s}\equiv
\left( x+_{\psi}B\right)^{n}\; ;\;\;n\geq 0.
\end{equation}
$\psi$- {\em basic Bernoulli-Ward polynomials} $\left\{
B_{n}(x)\right\}_{n\geq 0}$ are generalized Appell polynomials
i.e.
\begin{equation}\label{psi-B-3}
\partial_{\psi}B_{n}(x)=n_{\psi}B_{n-1}(x)
\end{equation}
and being $\psi$-Sheffer they naturally do satisfy the
$\psi$-Sheffer-Appell identity \cite{3,2}
\begin{equation}\label{psi-B-4}
B_{n}(x+_{\psi}y)=\sum_{s=0}^{n}\binom{n}{s}_{\psi}B_{s}(y)x^{n-s}.
\end{equation}
$\psi$- {\em basic Bernoulli-Ward polynomials} $\left\{
B_{n}(x)\right\}_{n\geq 0}$ are also equivalently characterized
via their $\psi$-exponential generating function
\begin{equation}\label{psi-B-5}
\sum_{n\geq
0}z^{n}\frac{B_{n}(x)}{n_{\psi}!}=\frac{z}{\exp_{\psi}\left\{
z\right\}-1}\exp_{\psi}\left\{ xz\right\}
\end{equation}
while the $\psi$-exponential generating function of
$\psi$-Bernoulli-Ward numbers $B_{n}=B_{n}(0) $ is
\begin{equation}\label{psi-B-6}
\sum_{n\geq
0}z^{n}\frac{B_{n}}{n_{\psi}!}=\frac{z}{\exp_{\psi}\left\{
z\right\}-1}\;\; ,B_{n}=B_{n}(0)\;\;\;n\geq 0.
\end{equation}
Compare the above with the theorem 16.2 in \cite{1}. There one
also shows that
\begin{equation}\label{psi-B-7}
\frac{B_{r+1}\left(
\overline{n}\right)-B_{r}}{(r+1)_{\psi}}=\sum_{k=0}^{n-1}\overline{k}^{r}\;\;
;r\geq 1.
\end{equation}
where
\begin{equation}\label{psi-B-8}
\overline{k}^{r}=\left( 1+_{\psi}1+_{\psi}\ldots
+_{\psi}1\right)^{r} \; \leftarrow k\;\; summands
\end{equation}
and $\psi$-{\em multinomial} formula reads (see \cite{1})
\begin{equation}\label{psi-B-9}
\left(x_{1}+_{\psi}x_{2}+_{\psi} \ldots +_{\psi}x_{k}\right)^{n}=
\sum_{\begin{array}{l} s_{1},\ldots s_{k}=0\\
s_{1}+s_{2}+\ldots +s_{k}=n
\end{array}}^{n}\binom{n}{s_{1},\ldots ,s_{k}}_{\psi}x_{1}^{s_{1}}\ldots x_{k}^{s_{k}}
\end{equation}
where
$$\binom{n}{s_{1},\ldots ,s_{k}}_{\psi}=\frac{n_{\psi}!}{(s_{1})_{\psi}!\ldots (s_{k})_{\psi}!}.$$
Naturally $\psi$- {\em basic Bernoulli-Ward polynomials} $\left\{
B_{n}(x)\right\}_{n\geq 0}$ satisfy the $\psi$-{\em difference}
equation
\begin{equation}\label{psi-B-10}
\Delta_{\psi}B_{n}(x)=n_{\psi}x^{n-1}\; ;n\geq 0
\end{equation}
hence they play the same role in $\psi$-{\em difference} calculus
as Bernoulli polynomials do in standard difference calculus (see:
Theorem 16.1 in \cite{1})due to the following:
 The central problem of the $Q(\partial_{\psi})$ -
{\em difference calculus} is:
$$ Q(\partial_{\psi})f=\varphi\;\;\;\;\;\;\;\;\;\;\varphi=?,$$
where $f, \varphi$ - are for example formal series of polynomials.

The idea of finding solutions is then following. As we know
\cite{2,3} any $\psi$- delta operator $Q$ is of the form
$Q(\partial_{\psi})=\partial_{\psi}\hat{B}$ where$\hat{B}\in
\Sigma_{\psi}$. Consider then $\Delta_{\psi}=\partial_{\psi}\hat{B
}\equiv
\partial_{\psi}\sum_{k\geq
0}\frac{1}{(k+1)_{\psi}}\frac{\partial_{\psi}^{k}}{k_{\psi}!}=E(\partial_{\psi})-id$,
hence ($\hat{B}\in\Sigma_{\psi}$) we have for $\hat{B}$-
recognized as $\psi$- Bernoulli operator - the obvious expression
$$\hat{B}=\frac{\partial_{\psi}}{\Delta_{\psi}}=\sum_{n\geq 0}\frac{B_{n}}{n_{\psi}!}\partial_{\psi}^{n}.$$
Now multiply by $\hat{B}\equiv
\frac{\partial_{\psi}}{e^{\partial_{\psi}}-1}\equiv \sum_{n\geq
0}\frac{B_{n}}{n_{\psi}!}\partial_{\psi}^{n}$ the equation
$\Delta_{\psi}f=\varphi$ in order to get
\begin{equation}\label{psi_B_11}
\partial_{\psi}f=\sum_{n\geq
0}\frac{B_{n}}{n_{\psi}!}\varphi^{(n)},\;\;\;\varphi^{(n)}=\partial_{\psi}\varphi^{(n-1)}.
\end{equation}
The solution then reads:
\begin{equation}\label{psi-B-12}
f(x)=\sum_{n \geq
1}\frac{B_{n}}{n_{\psi}!}\varphi^{(n-1)}(x)+\int_{\psi}\varphi(x)+p(x),
\end{equation}
where $p$ is  "$+_{\psi}1$- periodic" i.e. $p(x+_{\psi}1)=p(x)$
i.e. $\Delta_{\psi}p=0$. Here the $\psi$ - integration
$\int_{\psi}\varphi(x)$ is defined as in \cite{2}. We recall it in
brief. Let us introduce the following representation for
$\partial_{\psi}$ "difference-ization"
$$ \partial_{\psi}=\hat{n}_{\psi}\partial_{0}\; ;\;\;\;\hat{n}_{\psi}x^{n-1}=n_{\psi}x^{n-1};\;\;n\geq 1,$$
where $\partial_{0}x^{n}=x^{n-1}$ i.e.  $q=0$ "Jackson derivative"
$\partial_{0}$ is identical with divided difference operator. Then
we define the linear mapping $\int_{\psi}$ accordingly:
$$ \int_{\psi}x^{n}=\left( \hat{x}\frac{1}{\hat{n}_{\psi}}\right)x^{n}=\frac{1}{(n+1)_{\psi}}x^{n+1};\;\;\;n\geq 0$$
where of course  $\partial_{\psi}\circ \int_{\psi}=id$.
\section{Two Illustrative Specifications}
\subsection{$q$-umbral case \cite{1}-\cite{5}}
The following choice \cite{2,3,4,5} of the admissible sequence
$\psi_{n}(q)=\left[ R(q^{n})!\right]^{-1}$ and then
$R(x)=\frac{1-x}{1-q}$ results in the well known $q$-factorial
$n_{q}!=n_{q}(n-1)_{q}!,\;(n_{\psi}=n_{q})$ while the
$\psi$-derivative $\partial_{\psi}$ becomes the Jackson's
derivative $\partial_{q}$: $\left( \partial_{q}\varphi
\right)(x)=\frac{\varphi(x)-\varphi(qx)}{(1-q)x}$\cite{1}. The
$\psi$- integration \cite{2,5} becomes the well known $q$-
integration and we arrive at the $q$- Bernoulli numbers and $q$-
Bernoulli polynomials (for further references to Cigler, Roman and
others see\cite{2,3,4,5}).
\subsection{FFOC - case \cite{5}}
In straightforward analogy - (see FFOC-{\bf F}ibonomial {\bf
F}inite {\bf O}perator {\bf C}alculus, Example 2.1 in \cite{5}) -
consider now the Fibonomial coefficients ($F_{n}$ - Fibonacci
numbers)

$$\binom{n}{k}_{F}=\frac{F_{n}!}{F_{k}!F_{n-k}!}=\binom{n}{n-k}_{F},$$
where  $n_{F}\equiv F_{n}\neq 0$,
$n_{F}!=n_{F}(n-1)_{F}(n-2)_{F}(n-3)_{F}\ldots
2_{F}1_{F};\;\;0_{F}!=1$;
$$n_{F}^{\underline{k}}=n_{F}(n-1)_{F}\ldots
(n-k+1)_{F};\;\;\;\binom{n}{k}_{F}\equiv
\frac{n^{\underline{k}}_{F}}{k_{F}!}$$ and difference operator
$\partial_{F}$ linearly extended from
$\partial_{F}x^{n}=n_{F}x^{n-1};\;\;n\geq 0$ - we shall call the
$F$-derivative. Then in conformity with \cite{1} and with notation
as in \cite{2}-\cite{6} one writes:

\renewcommand{\labelenumi}{(\arabic{enumi})}
\begin{enumerate}
\item $\left( x+_{F}a\right)^{n}\equiv \sum_{k\geq
0}\binom{n}{k}_{F}a^{k}x^{n-k}$  where  $\binom{n}{k}_{F}\equiv
\frac{n^{\underline{k}}_{F}}{k_{F}!}$ \\and
$n_{F}^{\underline{k}}=n_{F}(n-1)_{F}\ldots (n-k+1)_{F}$; \item
$\left( x+_{F}a\right)^{n}\equiv
E^{a}(\partial_{F}x^{n};\;\;E^{a}(\partial_{F})=\sum_{n\geq
0}\frac{a^{n}}{n_{F}!}\partial_{F}^{n};\;\;\\E^{a}(\partial_{F})f(x)=f(x+_{F}a),$
$E^{a}(\partial_{F})$ is the corresponding generalized translation
operator.
\end{enumerate}

 The $\psi$- integration becomes now still
not explored $F$- integration and we arrive at the $F$- Bernoulli
numbers and $F$- Bernoulli polynomials - all to be investigated
soon.
\\
{\bf Note:} recently a combinatorial interpretation of Fibonomial
coefficient has been found \cite{6,7}.

 \end{document}